\documentclass[12pt]{article}
\usepackage{amssymb}
\usepackage{amsfonts}
\usepackage{latexsym}
\usepackage{amsthm}

\begin{document}


\renewcommand{\theequation}{\arabic{section}.\arabic{equation}}
\def\[{\begin{equation}}              \def\]{\end{equation}}
\def\lb{\label}                       \newcommand{\er}[1]{{\rm(\ref{#1})}}
\theoremstyle{plain}
\newtheorem{theorem}{\bf Theorem}[section]
\newtheorem{lemma}[theorem]{\bf Lemma}
\newtheorem{corollary}[theorem]{\bf Corollary}
\newtheorem{proposition}[theorem]{\bf Proposition}
\newtheorem{definition}[theorem]{\bf Definition}
\newtheorem{remark}[theorem]{\it Remark}

\def\a{\alpha}  \def\cA{{\cal A}}     \def\bA{{\bf A}}  \def\mA{{\mathscr A}}
\def\b{\beta}   \def\cB{{\cal B}}     \def\bB{{\bf B}}  \def\mB{{\mathscr B}}
\def\g{\gamma}  \def\cC{{\cal C}}     \def\bC{{\bf C}}  \def\mC{{\mathscr C}}
\def\G{\Gamma}  \def\cD{{\cal D}}     \def\bD{{\bf D}}  \def\mD{{\mathscr D}}
\def\d{\delta}  \def\cE{{\cal E}}     \def\bE{{\bf E}}  \def\mE{{\mathscr E}}
\def\D{\Delta}  \def\cF{{\cal F}}     \def\bF{{\bf F}}  \def\mF{{\mathscr F}}
\def\c{\chi}    \def\cG{{\cal G}}     \def\bG{{\bf G}}  \def\mG{{\mathscr G}}
\def\z{\zeta}   \def\cH{{\cal H}}     \def\bH{{\bf H}}  \def\mH{{\mathscr H}}
\def\e{\eta}    \def\cI{{\cal I}}     \def\bI{{\bf I}}  \def\mI{{\mathscr I}}
\def\p{\psi}    \def\cJ{{\cal J}}     \def\bJ{{\bf J}}  \def\mJ{{\mathscr J}}
\def\vT{\Theta} \def\cK{{\cal K}}     \def\bK{{\bf K}}  \def\mK{{\mathscr K}}
\def\k{\kappa}  \def\cL{{\cal L}}     \def\bL{{\bf L}}  \def\mL{{\mathscr L}}
\def\l{\lambda} \def\cM{{\cal M}}     \def\bM{{\bf M}}  \def\mM{{\mathscr M}}
\def\L{\Lambda} \def\cN{{\cal N}}     \def\bN{{\bf N}}  \def\mN{{\mathscr N}}
\def\m{\mu}     \def\cO{{\cal O}}     \def\bO{{\bf O}}  \def\mO{{\mathscr O}}
\def\n{\nu}     \def\cP{{\cal P}}     \def\bP{{\bf P}}  \def\mP{{\mathscr P}}
\def\r{\rho}    \def\cQ{{\cal Q}}     \def\bQ{{\bf Q}}  \def\mQ{{\mathscr Q}}
\def\s{\sigma}  \def\cR{{\cal R}}     \def\bR{{\bf R}}  \def\mR{{\mathscr R}}
\def\S{\Sigma}  \def\cS{{\cal S}}     \def\bS{{\bf S}}  \def\mS{{\mathscr S}}
\def\t{\tau}    \def\cT{{\cal T}}     \def\bT{{\bf T}}  \def\mT{{\mathscr T}}
\def\f{\phi}    \def\cU{{\cal U}}     \def\bU{{\bf U}}  \def\mU{{\mathscr U}}
\def\F{\Phi}    \def\cV{{\cal V}}     \def\bV{{\bf V}}  \def\mV{{\mathscr V}}
\def\P{\Psi}    \def\cW{{\cal W}}     \def\bW{{\bf W}}  \def\mW{{\mathscr W}}
\def\o{\omega}  \def\cX{{\cal X}}     \def\bX{{\bf X}}  \def\mX{{\mathscr X}}
\def\x{\xi}     \def\cY{{\cal Y}}     \def\bY{{\bf Y}}  \def\mY{{\mathscr Y}}
\def\X{\Xi}     \def\cZ{{\cal Z}}     \def\bZ{{\bf Z}}  \def\mZ{{\mathscr Z}}
\def\O{\Omega}

\def\ve{\varepsilon}   \def\vt{\vartheta}    \def\vp{\varphi}    \def\vk{\varkappa}

\def\Z{{\mathbb Z}}    \def\R{{\mathbb R}}   \def\C{{\mathbb C}}
\def\T{{\mathbb T}}    \def\N{{\mathbb N}}   \def\dD{{\mathbb D}}


\def\la{\leftarrow}              \def\ra{\rightarrow}            \def\Ra{\Rightarrow}
\def\ua{\uparrow}                \def\da{\downarrow}
\def\lra{\leftrightarrow}        \def\Lra{\Leftrightarrow}


\let\ge\geqslant                 \let\le\leqslant
\def\/{\over}                    \def\iy{\infty}
\def\sm{\setminus}               \def\es{\emptyset}
\def\ss{\subset}                 \def\ts{\times}
\def\pa{\partial}                \def\os{\oplus}
\def\om{\ominus}                 \def\ev{\equiv}
\def\iint{\int\!\!\!\int}        \def\iintt{\mathop{\int\!\!\int\!\!\dots\!\!\int}\limits}
\def\el2{\ell^{\,2}}             \def\1{1\!\!1}
\def\sh{\sharp}


\def\qq{\quad}                   \def\qqq{\qquad}
\def\lt{\biggl}                  \def\rt{\biggr}
\def\no{\noindent}               \def\ol{\overline}
\def\wt{\widetilde}              \def\wh{\widehat}
\newcommand{\nt}[1]{{\mathop{#1}\limits^{{}_{\,\bf{\sim}}}}\vphantom{#1}}
\newcommand{\nh}[1]{{\mathop{#1}\limits^{{}_{\,\bf{\wedge}}}}\vphantom{#1}}
\newcommand{\nc}[1]{{\mathop{#1}\limits^{{}_{\,\bf{\vee}}}}\vphantom{#1}}
\newcommand{\oo}[1]{{\mathop{#1}\limits^{\,\circ}}\vphantom{#1}}
\newcommand{\po}[1]{{\mathop{#1}\limits^{\phantom{\circ}}}\vphantom{#1}}


\def\Im{\mathop{\rm Im}\nolimits}
\def\Iso{\mathop{\rm Iso}\nolimits}
\def\Ker{\mathop{\rm Ker}\nolimits}
\def\Ran{\mathop{\rm Ran}\nolimits}
\def\Re{\mathop{\rm Re}\nolimits}
\def\Res{\mathop{\rm Res}\nolimits}
\def\Tr{\mathop{\rm Tr}\nolimits}
\def\arg{\mathop{\rm arg}\nolimits}
\def\const{\mathop{\rm const}\nolimits}
\def\det{\mathop{\rm det}\nolimits}
\def\diag{\mathop{\rm diag}\nolimits}
\def\dim{\mathop{\rm dim}\nolimits}
\def\dist{\mathop{\rm dist}\nolimits}
\def\rank{\mathop{\rm rank}\limits}
\def\res{\mathop{\rm res}\limits}
\def\sign{\mathop{\rm sign}\nolimits}
\def\span{\mathop{\rm span}\nolimits}
\def\supp{\mathop{\rm supp}\nolimits}

\newcommand\matr[4]{\left(\begin{array}{cc}#1 & #2 \cr #3 & #4\end{array}\right)}
\newcommand\vect[2]{\left(\begin{array}{c} #1 \cr #2 \end{array}\right)}

\newcommand\op[1]{\ol{1,#1}}


\title {Inverse spectral analysis for finite matrix-valued\\ Jacobi operators}

\author{Jochen Br\"uning\begin{footnote}{Institut f\"ur
Mathematik, Humboldt Universit\"at zu Berlin, e-mail: bruening@math.hu-berlin.de}
\end{footnote}
and Dmitry Chelkak\begin{footnote} {Dept. of Math. Analysis, Math. Mech. Faculty,
St.Petersburg State University. Universitetskij pr. 28, Staryj Petergof, 198504 St.Petersburg,
Russia, e-mail: delta4@math.spbu.ru }
\end{footnote}
and Evgeny Korotyaev\begin{footnote} {Correspondence author. Institut f\"ur Mathematik,
Humboldt Universit\"at zu Berlin, Rudower Chaussee 25, 12489 Berlin, Germany, e-mail:
evgeny@math.hu-berlin.de}
\end{footnote}}

\maketitle

\begin{abstract}
\no Consider the Jacobi operators $\cJ$ given by $(\cJ
y)_n=a_ny_{n+1}+b_ny_n+a_{n-1}^*y_{n-1}$, $y_n\in \C^m$ (here $y_0=y_{p+1}=0$), where
$b_n=b_n^*$ and $a_n:\det a_n\ne 0$ are the sequences of $m\ts m$ matrices, $n=1,..,p$. We
study two cases: (i) $a_n=a_n^*>0$; (ii) $a_n$ is a lower triangular matrix with real positive
entries on the diagonal (the matrix $\cJ$ is $(2m+1)$-band $mp\ts mp$ matrix with positive
entries on the first and the last diagonals). The spectrum of $\cJ$ is a finite sequence of
real eigenvalues $\l_1<\dots<\l_N$, where each eigenvalue $\l_j$ has multiplicity $k_j\le m$.
We show that the mapping $(a,b)\mapsto \{(\l_j,k_j)\}_1^N\oplus \{additional\ spectral\ data
\}$ is 1-to-1 and onto. In both cases (i), \nolinebreak (ii), we give the complete solution of
the inverse problem.
\end{abstract}

\section{Introduction  and main results}
\setcounter{equation}{0}

\no We consider the finite matrix-valued Jacobi operator $\cJ=\cJ_{a,b}$ acting in $(\C^p)^m$
and given by
\[
\label{Ja} \cJ=\left(\begin{array}{ccccccc} b_1 & a_1 & 0 & 0 &  ... &  0 \cr a_1^* & b_2 &
a_2 & 0 &  ... &  0 \cr 0 & a_2^* & b_3 & a_3 &  ... &  0 \cr ... & ... & ... & ... & ... &
...   \cr  0 & ...  & 0 & a_{p-2}^* & b_{p-1} & a_{p-1}\cr 0 & ...  & 0 & 0 & a_{p-1}^* & b_p
\end{array}\right),
\]
where $b=(b_n)_1^p$ and $a=(a_n)_1^{p-1}$ are the finite sequences of $m\ts m$ matrices such
that $b_n\in\bS=\{b:b=b^*\}$ and $\det a_n\ne 0$ for all $n\in \op{p}=\{1,2,\dots,p\}$. We
consider two cases:

\no (i)\,\ $a_n\in\bS_+= \{a:a=a^*>0\}$;

\no (ii) $a_n\in\bL_+=\{a: a$ is a lower triangular matrix with real positive entries on the
diagonal$\}$.

\pagebreak

In both cases we obtain the complete characterization of the set of spectral data that
correspond to these classes of operators $\cJ_{a,b}$. Note that in the second case $\cJ_{a,b}$
are self-adjoint $(2m\!+\!1)$-band matrices of the size $mp\ts mp$ with positive entries on
the first and the last diagonals.

The spectrum $\s(\cJ)$ of $\cJ=\cJ^*$ is a finite sequence of real eigenvalues
$$
\l_1<\l_2<\dots<\l_N,
$$
where each eigenvalue $\l_j$, $j\in \op N$ has multiplicity $k_j\in \op{m}$, i.e., $k_j$ is
the number of the linearly independent eigenvectors corresponding to $\l_j$. Note that
$$
k_1+k_2+\dots+k_N=mp.
$$
Introduce the fundamental $m\ts m$ matrix-valued solutions $\vp(z)=(\vp_n(z))_{0}^{p+1}$ and
$\c(z)=(\c_n(z))_{0}^{p+1}$ of the equation
\[
 \label{001}
 a_ny_{n+1}+b_ny_n+a_{n-1}^*y_{n-1}=z y_n,\qqq
 (z,n)\in\C\ts\op {p},
\]
such that
$$
\vp_{0}\ev \c_{p+1}\ev 0,\qq \vp_1\ev \c_{p}\ev I,
$$
where $I$ is the identity $m\ts m$ matrix and we set $a_0=a_p=I$ for convenience. Note that
$\vp_n(z)$ and $\c_n(z)$, $n\in\op{p}$ are matrix-valued polynomials such that
$\deg\vp_n(z)=n$, $\deg\c_n(z)=p\!+\!1\!-\!n$. Each eigenvector $\p=(\p_n)_1^p$,
$\cJ\p=\l_j\p$, has the form $\p_n=\vp_n(\l_j)v$ for some $v\in\C^m$. Hence, the eigenvalues
of $\cJ$ are zeros of $\det \vp_{p+1}(\cdot)$.

\begin{definition}[{\bf Spectral data}] For each eigenvalue $\l_j$, $j\in \op{N}$, we
define {the subspace}
\[
\label{cEDef} {\cE_j}=\Ker \vp_{p+1}(\l_j)= \lt\{h\in\C^m:\vp_{p+1}(\l_j)h=0\rt\}\ss\C^m,\qq
\dim\cE_j=k_j\le m,
\]
{\bf the orthogonal projector} $P_j:\C^N\to\cE_j\ss\C^N$ onto $\cE_j$ and {\bf the positive
self-adjoint operator} $g_j:\cE_j\to\cE_j$ given by
\[
\label{Gdef} {g_j}= G_j|_{\cE_j},\qq {where}\qq G_j=P_j\sum_{n=1}^{p} \vp_n^*(\l_j)
\vp_n(\l_j)P_j.
\]
\end{definition}

We now describe the connection between our spectral data and the matrix-valued Weyl-Titchmarsh
function $M(z)$ given by
\[
\label{mdef} M(z)=-\c_1(z)(\c_0(z))^{-1},\qqq z\in\C.
\]

\begin{proposition}
\label{Twt} The function $M(z)=M^*(\ol{z})$ is analytic in $\C\sm \{\l_j, j\in \op{N}\}$.
Moreover,
\[
\label{iwt} M(z)=-\sum_{j=1}^N {B_j\/z-\l_j},\qqq\qqq \sum_{j=1}^N B_j=I,
\]
where the self-adjoint matrices $B_j=\res_{z=\l_j}M(z)=B_j^*$ are given by
\[
\label{Bdef} B_j\big|_{\cE_j}=g_j^{-1},\qqq B_j\big|_{\C^N\om\cE_j}=0,\qq j\in\op{N}.
\]
\end{proposition}

In order to formulate our main result we need

\begin{definition}
We call the system $\{(\l_j,P_j),j\in\op{N}\}$ of the distinct real numbers $\l_j$ and the
orthogonal projectors $P_j:\C^m\to\C^m$ {\bf the $\bf p$\,--\,tame system}, if
$\sum_{j=1}^N\rank P_j=mp$ and
$$
\det \left(\begin{array}{cccc} T_0 & T_1 & ... & T_{p-1} \cr T_1 & T_2 & ... & T_p \cr ... &
... & ... & ... \cr T_{p-1} & T_{p} & ... & T_{2p-2}\end{array}\right) \ne 0,\qqq
T_s=\sum_{j=1}^N\l_j^sP_j,\qq s=0,..,{2p\!-\!2}.
$$
\end{definition}

\noindent {\it Remark. } (i) The $mp\ts mp$ matrix $(T_{sk})_{s,k=0}^{p-1}$ is always
non-negative definite (see (\ref{Tge0})). Hence, the system $\{(\l_j,P_j),j\in\op{N}\}$ is
$p$\,--tame iff this matrix is strictly positive definite.

 \no (ii) Let $N=p$ and $k_j=\dim \cE_j = \rank P_j = m$ for
all $j\in\op{p}$, i.e., $P_j=I$ for all $j$. Then for each distinct values $\l_1<\dots<\l_p$
the system $\{(\l_j,I),j\in\op{p}\}$ is $p$\,--\,tame.

\no (iii) Let $p=2$, $N=3$, $k_1+k_2=k_3=m$. Then the system
$\{(\l_1,P_1),(\l_2,P_2),(\l_3,I)\}$ is $2$\,--\,tame iff $\Ker P_1\cap\Ker P_2=\{0\}$.

\medskip

Introduce the set of spectral data
\[
\label{Sdef} \cS_p= \left\{ \left(\l_j,P_j,g_j\right)_1^N: \begin{array}{l} p\le N\le pm,\
\l_1<\dots<\l_N {\rm \ are\ real\ numbers}; \cr \{(\l_j,P_j),j\in\op{N}\}\ {\rm is\ a\ }
p{\rm\,-\,tame\ system}; \cr g_j:\Ran P_j\to\Ran P_j\ {\rm are\ linear\ operators\ such} \cr
{\rm that\ } g_j=g_j^*>0\ {\rm and\ } \sum_{j=1}^N P_jg_j^{-1}P_j = I
\end{array}\right\}
\]
and the mapping
$$
\Psi: \left((a_n)_1^{p-1};(b_n)_1^p\right)\mapsto \left(\l_j,P_j,g_j\right)_1^N\,.
$$

We formulate our main result.
\begin{theorem}\label{Thm}
(i) The mapping $\Psi:\bS_+^{p-1}\!\ts\bS^p\to\cS_p$ is one-to-one and onto.

\no (ii) The mapping $\Psi:\bL_+^{p-1}\!\ts\bS^p\to\cS_p$ is one-to-one and onto.
\end{theorem}

There is an enormous literature on inverse spectral problems for scalar (i.e., $m=1$) Jacobi
matrices (see book \cite{T} and references therein), but considerably less for matrix-valued
Jacobi operators (see \cite{CGR} and references therein). The inverse problems for finite
scalar Jacobi matrices were considered by several authors (see \cite{GS} and references
therein). Some uniqueness results for matrix-valued Jacobi operators were obtained in
\cite{CGR}, and the intimate connection to matrix-valued orthogonal polynomials and the moment
problem were treated in \cite{L},\cite{DL}.

The main goal of our paper is to give the complete characterization of the set of spectral
data for operators $\cJ_{a,b}$. We hope to use similar ideas to obtain the complete
characterization of the spectral data for the Sturm-Liouville operators $\cH y=-y''+V(x)y$,
$y(0)=y(1)=0$, on the unit interval $[0,1]$, where $V=V^*$ is a $m\ts m$ matrix potential.
Note that the "local characterization" of the spectral data for operators $\cH$ was obtained
in \cite{CK}.

In the proof we use the approach from \cite{GS} and some technique from \cite{CK}. The inverse
spectral problem consists of the following parts: i) \nolinebreak uniqueness, ii) \nolinebreak
characterization, iii) \nolinebreak reconstruction. In Theorem \ref{Thm} we solve all these
problems simultaneously.

\section{Preliminaries}
\setcounter{equation}{0}

For each two sequences of polynomials $\vt(z)=(\vt_n(z))_0^{p+1}$, $\e=(\e_n(z))_0^{p+1}$ we
introduce the Wronskian
$$
\{\vt,\e\}_n(z)=\vt_n^*(\ol{z}) a_n \e_{n+1}(z)- \vt_{n+1}^*(\ol{z}) a^*_n\e_{n}(z),\qqq
n=0,..,p.
$$
If (\ref{001}) holds for both $\vt$ and $\e$, then $\{\vt,\e\}_n(z)$ does not depend on $n$.
In particular,
\[
\label{wron} \vp_p^*(\ol{z})\vp_{p+1}(z)-\vp_{p+1}^*(\ol{z})\vp_p(z)=\{\vp,\vp\}(z)=0\qq {\rm
and}\qq \c_0^*(\ol{z})= \{\c,\vp\}(z)=\vp_{p+1}(z),
\]
since $a_0=a_p=I$ by our convention. Recall that
$$
\cE_j=\Ker\vp_{p+1}(\l_j)=\{h\in \C^m: \vp_{p+1}(\l_j)h=0\},\qq j\in\op{N},
$$
and $P_j:\C^m\to\cE_j$ is the orthogonal projector. Also, we need the subspaces
$$
\cE_j^\sh=\Ker\vp_{p+1}^*(\l_j)=\{h\in \C^m: \vp_{p+1}^*(\l_j)h=0\},\qq j\in\op{N},
$$
and the orthogonal projectors $P_j^\sh:\C^m\to\cE_j^\sh$. Note that
$\dim\cE_j^\sh=\dim\cE_j=k_j$.

\begin{lemma} For each $j\in\op{N}$ the following identities are fulfilled:
\[
\label{i1} P_j^\sh\vp_p(\l_j)P_j=\vp_p(\l_j)P_j,\qqq \c_1(\l_j)\vp_p(\l_j)P_j=P_j,
\]
\[
\label{s2} I\le \sum_{j=1}^p \vp_n^*(\l_j)\vp_n(\l_j) =
\vp_p^*(\l_j)\dot{\vp}_{p+1}(\l_j)- \vp_{p+1}^*(\l_j)\dot{\vp}_p(\l_j).
\]
\end{lemma}
\begin{proof}
Using (\ref{wron}), we obtain
$\vp_{p+1}^*(\l_j)\vp_p(\l_j)P_j=\vp_p^*(\l_j)\vp_{p+1}(\l_j)P_j=0$. This means
$\vp_p(\l_j)P_j\ss\cE_j^\sh$, i.e., $P_j^\sh\vp_p(\l_j)P_j=\vp_p(\l_j)P_j$. Let
$$
y_n=\c_n(\l_j)\vp_p(\l_j) P_j-\vp_n(\l_j)P_j,\qq n=0,..,p+1.
$$
The sequence $(y_n)_0^{p+1}$ satisfies (\ref{001}) for $z=\l_j$ and
$$
y_{p+1}= - \vp_{p+1}(\l_j)P_j=0,\qqq y_{p}=\vp_p(\l_j)P_j - \vp_p(\l_j)P_j=0.
$$
This yields $y_n=0$ for all $n\in\op{p}$. In particular, $y_1=\c_1(\l_j)\vp_p(\l_j)
P_j-P_j=0$.

Furthermore, using equation (\ref{001}), we obtain
$$
\vp_{n+1}^*(\ol{z}) a_n^* + \vp_n^*(\ol{z})(b_n-z) + \vp_{n-1}^*(\ol{z})a_{n-1}=0,
$$
$$
a_n\dot\vp_{n+1}(z) + (b_n-z)\dot\vp_n(z) + a_{n-1}^*\dot\vp_{n-1}(z)=\vp_n(z).
$$
Multiplying the first equation by $\dot{\vp}_n(z)$ from the right, the second equation by
$\vp_n^*(\ol{z})$ from the left and taking the difference, we deduce that
$$
\vp_n^*(\ol{z})\vp_n(z)=
\{\vp, \dot \vp\}_{n}(z) -\{\vp, \dot \vp\}_{n-1} (z),\qq n\in\op{p}.
$$
The summing implies (\ref{s2}) since $\{\vp,\dot{\vp}\}_0=0$. Note that $I=
\vp_1^*(\l_j)\vp_1(\l_j)$ by definition.
\end{proof}

\begin{lemma}
\label{L^-1} (i) For each $j\in\op{N}$ the mapping
$Y_j=P_j^\sh\dot{\vp}_{p+1}(\l_j)P_j:\cE_j\to\cE_j^\sh$ is invertible.

\no (ii) For each $j\in\op{N}$ the following asymptotics holds:
\[
\label{VpRes} (\vp_{p+1}(z))^{-1} = \frac{P_jY_j^{-1}P_j^\sh}{z-\l_j}+O(1)\qq {as}\qq
z\to\l_j.
\]
\end{lemma}
\begin{proof}
(i) Let $Y_j h =0$ for some $h\in \cE_j$. Using (\ref{s2}) and (\ref{i1}), we obtain
$$
P_j\le P_j \sum_{n=1}^p \vp_n^*(\l_j)\vp_n(\l_j) P_j =
P_j\vp_p^*(\l_j)\dot{\vp}_{p+1}(\l_j)P_j = P_j\vp_p^*(\l_j)P_j^\sh\cdot
P_j^\sh\dot{\vp}_{p+1}(\l_j)P_j.
$$
Since the left hand side is positive definite on $\cE_j$, the operator
$Y_j=P_j^\sh\dot{\vp}_{p+1}(\l_j)P_j$ is invertible.

\no (ii) Let
$$
\vp_{p+1}(z) = \matr {A(z)} {B(z)} {C(z)} {D(z)},
$$
where
$$
\begin{array}{llll}
A=(P_j^\sh)^\bot\vp_{p+1}P_j^\bot:&\cE_j^\bot\to(\cE_j^\sh)^\bot,\qq &
B=(P_j^\sh)^\bot\vp_{p+1}P_j:&\cE_j\to(\cE_j^\sh)^\bot, \cr
C=P_j^\sh\vp_{p+1}P_j^\bot:&\cE_j^\bot\to\cE_j^\sh\ , &
D=P_j^\sh\vp_{p+1}P_j:&\cE_j\to\cE_j^\sh\ .
\end{array}
$$
Since $\vp_{p+1}(\l_j)P_j=0$ and $P_j^\sh\vp_{p+1}(\l_j)=0$, we have $B,C=O(t)$ as
$t=z-\l_j\to 0$. Moreover,
$$
D = t\cdot P_j^\sh\dot{\vp}_{p+1}(\l_j)P_j+O(t^2)= tY_j + O(t^2),\qqq A = X_j+O(t),
$$
where $X_j=(P_j^\sh)^\bot\vp_{p+1}(\l_j)P_j^\bot=\vp_{p+1}(\l_j)P_j^\bot:
\cE_j^\bot\to(\cE_j^\sh)^\bot$. The operator $X_j$ is invertible, since $\vp_{p+1} h=0$
implies $h \in \cE_j$. The operator $Y_j$ is invertible due to (i). Therefore, the standard
formula for the inverse matrix
$$
{\matr A B C D}^{-1}=\matr {A^{-1}+A^{-1}BH^{-1}CA^{-1}} {-A^{-1}BH^{-1}} {-H^{-1}CA^{-1}}
{H^{-1}},\qq H=D-CA^{-1}B,
$$
gives
$$
(\vp_{p+1}(z))^{-1} = \matr {X_j^{-1}+O(t)} {O(1)} {O(1)} {t^{-1}Y_j^{-1}+O(1)} \qq {\rm as}
\qq t\to 0.
$$
In particular, $\res_{z=\l_j}(\vp_{p+1}(z))^{-1}=P_jY_j^{-1}P_j^\sh$.
\end{proof}

\begin{proof}[{\bf Proof of Proposition \ref{Twt}.}]
Note that $\c_0^*(\ol{z})\c_1(z)-\c_1^*(\ol{z})\c_0(z)=\{\c,\c\}_0(z)=0$. This gives
$M^*(\ol{z})=M(z)$, $z\in\C$. Due to (\ref{VpRes}), the function
$(\c_0(z))^{-1}=(\vp_{p+1}^*(\ol{z}))^{-1}$ has a simple pole at each point $z=\l_j$.
Therefore, the function $M(z)=-\c_1(z)(\c_0(z))^{-1}$ has a simple pole at $z=\l_j$ and
$$
B_j=-\res_{z=\l_j}M(z)=\c_1(\l_j)\cdot (\res_{z=\l_j} (\vp_{p+1}(z))^{-1})^* = \c_1(\l_j)
P_j^\sh (Y_j^*)^{-1} P_j,
$$
where $Y_j=P_j^\sh\dot{\vp}_{p+1}(\l_j)P_j:\cE_j\to\cE_j^\sh$. This yields
$B_j\big|_{\C^N\om\cE_j}=0$. Furthermore, using (\ref{s2}) and (\ref{i1}), we obtain
$$
G_j= G_j^* = P_j\dot \vp_{p+1}^*(\l_j)\vp_p(\l_j) P_j = P_j\dot\vp_{p+1}^*(\l_j) P_j^\sh
\vp_p(\l_j)P_j = P_j Y_j^* P_j^\sh \cdot \vp_p(\l_j)P_j.
$$
Hence (\ref{i1}) yields,
$$
B_j G_j = \c_1(\l_j)P_j^\sh\vp_p(\l_j)P_j = \c_1(\l_j)\vp_p(\l_j)P_j = P_j,
$$
i.e., $B_j\big|_{\cE_j}=g_j^{-1}$. Asymptotics $M(z)=-Iz^{-1}+O(z^{-2})$ as $z\to\iy$ (see
(\ref{001}) for $n=0$) and the standard arguments from the function theory give (\ref{iwt}).
\end{proof}

\section{Proof of Theorem \ref{Thm}}
\setcounter{equation}{0}

\begin{lemma}
\label{Mom} The system $\{(\l_j,P_j),j\in\op{N}\}$ is $p$\,--\,tame, iff there is no nonzero
vector-valued polynomial $F(z)\in\C^m$ such that $\deg F\le p\!-\!1$ and $P_jF(\l_j)=0$ for
all $j\in \op{N}$.
\end{lemma}

\begin{proof}
Note that for each vector $v=(v_s)_0^{p-1}$, $v_s\in\C^m$, we have
$$
\left(v_0^*,v_1^*,...,v_{p-1}^*\right) \left(\begin{array}{cccc} T_0 & T_1 & ... & T_{p-1} \cr
T_1 & T_2 & ... & T_p \cr ... & ... & ... & ... \cr T_{p-1} & T_{p} & ... &
T_{2p-2}\end{array}\right) \left(\begin{array}{c} v_0\cr v_1\cr ... \cr
v_{p-1}\end{array}\right)
$$
\[
\label{Tge0} =\sum_{s,k=0}^{p-1} v_s^*\sum_{j=1}^N \l_j^{s+k} P_j v_k = \sum_{j=1}^N
F^*(\l_j)P_j F(\l_j)\ge 0,
\]
where $F(z)=\sum_{s=0}^{p-1}z^sv_s$. Therefore, $\det (T_{s+k})_{s,k=0}^{p-1} \ne 0$ iff there
is no nonzero polynomial $F(z)$ of degree at most $p\!-\!1$ such that $P_jF(\l_j)=0$ for all
$j\in\op{N}$.
\end{proof}

\begin{proof}[{\bf Proof of Theorem \ref{Thm}.}]
The proof is similar for both cases (i),(ii) and consists of three parts. We need also two
simple technical Lemmas \ref{Hergl}, \ref{Sqrt} that are located at the end.

\medskip

\no {\bf 1. $\bf\Psi$ maps $\bf \bS_+^{p-1}\!\ts\bS^p$ and $\bf \bL_+^{p-1}\!\ts\bS^p$ into
$\bf\cS_p$.}

\no Recall that the identity $\sum_{j=1}^N P_jg_j^{-1}P_j = \sum_{j=1}^N B_j = I$ is proved in
Proposition \ref{Twt}. In order to check that $\{(\l_j,P_j),j\in\op{N}\}$ is $p$\,--\,tame, we
use Lemma \ref{Mom}. Suppose that $P_jF(\l_j)=0$ for all $j\in \op{N}$ and some vector-valued
polynomial $F(z)$ of degree at most $p-1$. Using asymptotics (\ref{VpRes}) near poles of
$(\c_0(z))^{-1}=(\vp_{p+1}^*(\ol{z}))^{-1}$, we deduce that the vector-valued function
$(\c_0(z))^{-1}F(z)$ is entire. On the other hand, $\c_0^{-1}(z)F(z)=O(z^{-p}\cdot
z^{p-1})=O(z^{-1})$ as $z\to\iy$. The Liouville Theorem gives $\c_0^{-1}(z)F(z)=0$ for all
$z\in\C$ and hence $F=0$.

\medskip

\no {\bf 2. $\bf\Psi:\bS_+^{p-1}\!\ts\bS^p\to\cS_p$ is one-to-one and
$\bf\Psi:\bL_+^{p-1}\!\ts\bS^p\to\cS_p$ is one-to-one.}

\no Following \cite{GS}, we introduce the sequence of $M$-functions
$$
M_n(z)= -\c_{n}(z)[a_{n-1}^*\c_{n-1}(z)]^{-1},\qq n\in\op{p\!+\!1}.
$$
Using (\ref{001}), it is easy to see that
\[
\label{MIdent} -(M_n(z))^{-1}=[(z-b_n)\c_n(z)-a_n\c_{n+1}(z)](\c_n(z))^{-1}= Iz-b_n + a_n
M_{n+1}(z) a_n^*
\]
for all $n\in\op{p}$. In particular, $M_n(z)=-z^{-1}+O(z^{-2})$ as $z\to\infty$ for all $n$.
Due to Proposition \nolinebreak \ref{Twt}, the spectral data $(\l_j,P_j,g_j)_1^N$ uniquely
determine the function $M=M_0$. Therefore, the matrices $b_1$ and $A=a_1a_1^*$ are uniquely
determined by the asymptotics
$$
M(z)=-Iz^{-1}-b_1z^{-1}-(a_1a_1^*+b_1^2)z^{-2}+O(z^{-3})\qq {\rm as}\qq z\to\iy.
$$
If $a_1,\wt{a}_1\in\bS_+$ (the case (i)) and $a_1a_1^*=A=\wt{a}_1\wt{a}_1^*$, then
$a_1=\wt{a}_1$. The same is true, if $a_1,\wt{a}_1\in\bL_+$ (the case (ii), see Lemma
\ref{Sqrt}). Hence, the matrices $b_1$ and $a_1$ are uniquely determined by the spectral data.
Therefore, the function $M_1(z)=a_1^{-1}(Iz-b_1-(M(z))^{-1})(a_1^{-1})^*$ is uniquely
determined too. Repeating this procedure, one uniquely determines all matrices $b_n$, $a_n$ in
both cases (i), (ii).

\medskip

\no {\bf 3. $\bf\Psi:\bS_+^{p-1}\!\ts\bS^p\to\cS_p$ is onto and
$\bf\Psi:\bL_+^{p-1}\!\ts\bS^p\to\cS_p$ is onto.}

\no Let $(\l_j,P_j,g_j)_1^N\in\cS_p$. We shall construct $b_n$ and $a_n$ step by step.
Introduce the function
\[
\label{M=} M(z)= -\sum_{j=1}^N{B_j\/z-\l_j},
\]
where $B_j=P_jg_j^{-1}P_j=B_j^*\ge 0$. Note that $M(\ol{z})=M^*(z)$, $z\in\C$, and
$$
\Im M(z)= \frac{1}{2i}(M(z)-M^*(z))>0\qq {\rm for\ all}\qq \Im z>0.
$$
In particular, $\det M(z)\ne 0$, if $\Im z>0$. This gives
$$
-\Im M^{-1}(z) = (M^{-1})^*(z)\cdot \Im M(z)\cdot M^{-1}(z)>0,\qq \Im z>0.
$$
Moreover, $-M^{-1}(z)=Iz+O(1)$ as $z\to\iy$, since $\sum_{j=1}^N B_j = I$. Therefore, the
Herglotz representation theorem for rational matrix-valued functions (see Lemma \ref{Hergl})
yields
\[
\label{M^-1=} -M^{-1}(z) = Iz + C - \sum_{s=1}^{K} {D_s\/z-\m_s}
\]
for some $\m_s\in \R$, $C=C^*$ and $D_s=D^*_s\ge 0$, $s\in \op{K}$. Denote $\cF_s=\C^N\om \Ker
D_s$ and let $Q_s:\C^N\to\cF_s\ss\C^N$ be the orthogonal projectors onto $\cF_s$. In
accordance with (\ref{MIdent}), we set
$$
b_1=-C.
$$

Let $\l_j\ne\m_s$ for all $j,s$. Then all poles of the meromorphic function $\det M(z)$ are
$\{\l_j\}_1^N$ and all roots are $\{\m_s\}_1^{K}$. Moreover, each $\l_j$ is a pole of the
multiplicity $\dim\cE_j=\rank P_j$ and each $\m_s$ is a root of the multiplicity
$\dim\cF_s=\rank Q_s$. Since $\det M(z)=(-z)^{-m}+O(z^{-m-1})$ as $z\to \iy$, we have
\[
\label{pr} \sum_{s=1}^{K}\rank Q_s = -m+\sum_{j=1}^N \rank P_j  = m(p-1).
\]
If $\l_j=\m_s$ for some $j,s$, then the corresponding pole and the root (partially) compensate
each other but (\ref{pr}) still holds true.

If $p=1$, then (\ref{pr}) yields $K=0$ and the reconstruction procedure stops.

Let $p>1$. We show that the system $\{(\m_s,Q_s),s\in\op{K}\}$ is $(p-1)$--\,tame, this is the
crucial point of the proof. Suppose that $Q_sG(\m_s)=0$ for all $s$ and some vector-valued
polynomial $G(z)$ of degree at most $p\!-\!2$. Due to Lemma \ref{Mom}, it is sufficient to
prove that $G=0$. In view of (\ref{M^-1=}), the vector-valued function
$$
F(z)=M^{-1}(z)G(z)
$$
is entire and $F(z)=O(z\cdot z^{p-2})=O(z^{p-1})$ as $z\to\infty$, i.e., $F(z)$ is a
vector-valued polynomial of degree at most $p\!-\!1$. Since $M(z)F(z)=G(z)$ is entire too,
(\ref{M=}) yields $P_jF(\l_j)=0$ for all $j\in\op{N}$. Due to Lemma \ref{Mom}, this implies
$F=0$. Hence, $G=0$.

Define $A=\sum_{s=1}^{K}D_s$ (note that $A=A^*\ge 0$). In order to show that $A>0$, suppose
that $f^* A f = 0$ for some constant vector $f\in\C^m$, $f\ne 0$. This gives $Q_s f=0$ for all
$s\in\op{K}$, which is a contradiction with Lemma \ref{Mom}. Hence, $A>0$. In accordance with
(\ref{MIdent}), we choose $a_1$ such that
$$
a_1a_1^*=A.
$$
In both cases: (i) $a_1=\sqrt{A}>0$ and (ii) see Lemma \ref{Sqrt}, there exists a unique
matrix solution $a_1$ of this equation. Now we set
$$
M_1(z) = -a_1^{-1}(b_1-Iz-(M(z))^{-1})(a_1^{-1})^*= - \sum_{s=1}^{K} {\wt{D}_s\/z-\m_s}\,,
$$
where
$$
\wt{D}_s=a_1^{-1}D_s(a_1^{-1})^*,\qq s\in\op{K},\qqq {\rm and}\qqq \sum_{s=1}^K \wt{D}_s=I.
$$

Note that $\Ker \wt{D}_s = a_1^*\Ker D_s=a_1^*\cF_s^\bot$. Let
$\wt{Q}_s:\C^N\to(a_1^*\cF_s^\bot)^\bot$ be the orthogonal projectors. Suppose that $\wt{Q}_s
G(\m_s)=0$ for all $s\in\op{K}$ and some vector-valued polynomial $G(z)$ of degree at most
$p\!-\!2$. Then $G(\m_s)\in a_1^*\cF_s^\bot$, i.e., $Q_s[(a_1^*)^{-1}G(\m_s)]=0$ for all
$s\in\op{K}$. Due to Lemma \ref{Mom}, this gives $(a_1^*)^{-1}G=0$ and hence $G=0$. Therefore,
the new system $\{(\m_s,\wt{Q}_s),s\in\op{K}\}$ is $(p\!-\!1)$--\,tame. Repeating the
procedure given above, we reconstruct $b_2$, $a_2$, $b_3$, $a_3$ and so on.
\end{proof}

\smallskip

\begin{lemma}
\label{Hergl} Let the rational matrix-valued function $f(z)$ have only real poles and
satisfies $f(z)=f^*(\ol{z})$, $z\in\C$, and $\Im f(z) = {1\/2i}\,(f(z)\!-\!f^*(z))>0$, $\Im
z>0$. Then $f(z)$ has the form
$$
f(z)=D_0 z + C - \sum_{s=1}^K {D_s\/z-\l_s}
$$
for some $K\ge 0$, $\l_s\in\R$, $C=C^*$ and $D_s=D_s^*\ge 0$, $s=0,..,K$.
\end{lemma}
\begin{proof}
Due to the identity $f(z)=f^*(\ol{z})$, we have
$$
f(z)= \sum_{j=1}^{K_0}D_{0,j} z^{j+1} + C + \sum_{s=1}^K \sum_{j=1}^{K_j}
{D_{s,j}\/(z-\l_s)^j}
$$
for some $D_{s,j}=D_{s,j}^*$ and $C=C^*$. Using the condition $\Im f(z)>0$ as $\Im z>0$ near
the points $z=\l_s$ and $z=\infty$, we deduce that $D_{s,j}=0$, if $j \ge 2$, and
$D_s=-D_{s,1}\ge 0$.
\end{proof}

\begin{lemma}
\label{Sqrt} Let $A\in\bS_+$. Then there exists unique matrix $a\in\bL_+$ such that $aa^*=A$.
\end{lemma}

\begin{proof}
Let
$$
A=\matr {B} {C^*} {C} {D},\qqq a=\matr {b} {0} {c} {d}
$$
where $B$,$b$ are positive real numbers, $C$,$c$ are vectors and $D$,$d$ are $(m-1)\ts(m-1)$
matrices. The equation $aa^*=A$ is equivalent to
$$
b\ol{b}=b^2=B,\qqq c\ol{b}=C,\qqq dd^*=D-cc^*=D-{CC^*\/B}.
$$
Note that $D-{1\/B}\,CC^*>0$, since $A>0$. Now the problem is reduced to the similar problem
for $(m-1)\ts(m-1)$ matrices. The case $m=1$ is trivial.
\end{proof}

\end{document}